\documentclass{amsproc}
\usepackage{amssymb}
\usepackage{graphicx}

  \def\NN{{\mathbb N}}

  \def\ZZ{{\mathbb Z}}
  \def\RR{{\mathbb R}}
  \def\kk{{\mathbf k}}

  \def\dD{{\mathcal D}}

  \def\sS{{\mathcal S}}
  \def\wW{{\mathcal W}}
  \def\csum{{\rm csum}}
  
  \def\asc{{\rm asc}}

  \def\ex{{\rm exc}}
  
  \def\fex{{\rm fexc}}
  
  \def\link{{\rm link}}

  \def\sd{{\rm sd}}

  \def\sm{\smallsetminus}
  \def\qed{$\hfill \Box$}

\newtheorem{theorem}{Theorem}[section]
\newtheorem{conjecture}[theorem]{Conjecture}
\newtheorem{corollary}[theorem]{Corollary}
\newtheorem{problem}[theorem]{Problem}
\newtheorem{question}[theorem]{Question}

\theoremstyle{definition}
\newtheorem{definition}[theorem]{Definition}
\newtheorem{example}[theorem]{Example}

\theoremstyle{remark}
\newtheorem{remark}[theorem]{Remark}

\numberwithin{equation}{section}

\begin{document}

\title{A survey of subdivisions and local $h$-vectors}

\author{Christos~A.~Athanasiadis}
\address{Department of Mathematics, University of Athens, Athens 15784,
Hellas (Greece)}
\email{caath@math.uoa.gr}

%    General info
\date{March 25, 2014 and, in revised form, October 20, 2015.}

\dedicatory{Dedicated to Richard Stanley on the occasion of his 70th
birthday.}


\begin{abstract}
The enumerative theory of simplicial subdivisions (triangulations) of
simplicial complexes was developed by Stanley in order to understand
the effect of such subdivisions on the $h$-vector of a simplicial complex.
A key role there is played by the concept of a local $h$-vector. This
paper surveys some of the highlights of this theory and recent
developments, concerning subdivisions of flag homology spheres and
their $\gamma$-vectors. Several interesting examples and open problems
are discussed.
\end{abstract}

\maketitle

\section{Introduction}
\label{sec:intro}

Let $\Delta$ be a $(d-1)$-dimensional simplicial complex (all such
complexes considered in this paper will be assumed to be finite). Much of the enumerative combinatorics of $\Delta$ is captured in the sequence $f(\Delta) = (f_{-1}(\Delta), f_0(\Delta),\dots,f_{d-1}(\Delta))$, called the \emph{$f$-vector} of $\Delta$, where $f_i (\Delta)$ denotes the number of $i$-dimensional faces of $\Delta$.
It is often more convenient to work with the \emph{$h$-vector} $h(\Delta) = (h_0(\Delta), h_1(\Delta),\dots,h_d(\Delta))$ of
$\Delta$, defined by
\begin{equation} \label{eq:hdef}
h(\Delta, x) \ = \ \sum_{i=0}^d h_i(\Delta) \, x^i \ = \
\sum_{i=0}^d f_{i-1} (\Delta) \, x^i (1-x)^{d-i},
\end{equation}
rather than directly with $f(\Delta)$, as several fundamental properties of the latter can be expressed in terms of $h(\Delta)$ in
an elegant way. For instance, the coordinates of $h(\Delta)$ are nonnegative for every Cohen-Macaulay complex $\Delta$ (in particular, whenever $\Delta$ is a homology ball or sphere) and symmetric for
every homology sphere $\Delta$; we refer to~\cite{StaCCA} for an authoritative exposition, where important algebraic and combinatorial interpretations of $h(\Delta)$ can also be found. When $\Delta$ has a nice combinatorial description, $h(\Delta)$ often turns out to be an enumerative invariant which is simpler and more basic than
$f(\Delta)$; typical situations appear in Section~\ref{sec:ex}. The
polynomial which appears in (\ref{eq:hdef}) is known as the
\emph{$h$-polynomial} of $\Delta$.

Motivated in part by analogous questions in Ehrhart theory, Kalai
and Stanley (see \cite{Sta92}) were interested to understand how the
$h$-vector of a simplicial complex changes after arbitrary simplicial subdivision. Specifically, they raised the question whether
$h(\Delta)$ increases coordinatewise after such subdivision of a
Cohen-Macaulay complex $\Delta$ (this monotonicity is obvious for the sum of the coordinates of $h(\Delta)$, which equals $f_{d-1}(\Delta)$, but highly nontrivial for each $h_i (\Delta)$). This paper is a brief survey of the resulting theory \cite{Sta92}, in which the concept of a local $h$-vector plays a prominent role, as well as of some recent developments. Our goal is neither to cover this admirable work of Richard Stanley in full generality, nor to convince the reader of
its importance and depth, but rather to convey part of its beauty through the mysterious and attractive properties of its main characters, the wealth of intriguing open problems it hides and its relevance to some popular topics in enumerative and geometric combinatorics.

After defining the various notions of simplicial subdivision which are relevant in the theory, Section~\ref{sec:sub} introduces local
$h$-vectors, lists their main properties and explains Stanley's beautiful solution to the monotonicity problem for $h(\Delta)$. The concept of a relative local $h$-vector, introduced in \cite{Ath12a}, is also discussed there. Section~\ref{sec:flag} discusses extensions and applications \cite{Ath12a} of the theory to subdivisions of flag homology spheres and the analogous monotonicity problem for a finer enumerative invariant than the $h$-vector, introduced in \cite{Ga05}, namely the $\gamma$-vector. Section~\ref{sec:ex} includes examples for which the computation of the local $h$-vector leads to interesting problems on the combinatorics of objects such as words, colored permutations and noncrossing partitions. Several open problems are included.

Due to lack of space, neither the general theory of subdivisions of lower graded posets developed in \cite[Part II]{Sta92}, nor its analogues to cubical $h$-vectors of cubical complexes \cite{Ath12b}, will be treated here.

\section{Simplicial subdivisions and local $h$-vectors}
\label{sec:sub}

This section summarizes the enumerative theory of simplicial subdivisions of \cite[Part I]{Sta92}. Some examples and concepts which first appeared in \cite{Ath12a} are also discussed. Background on simplicial complexes can be found in \cite{Bj95}.

\subsection{Subdivisions} \label{subsec:sub}
Consider two geometric simplicial complexes $\Sigma'$ and $\Sigma$, with corresponding abstract simplicial complexes $\Delta'$ and
$\Delta$, such that $\Sigma'$ subdivides $\Sigma$ (this means that every simplex of $\Sigma'$ is contained in some simplex of $\Sigma$ and that the union of the simplices of $\Sigma'$ is equal to the union of the simplices of $\Sigma$). Then, given any simplex $L' \in \Sigma'$, the relative interior of $L'$ is contained in the relative interior of a unique simplex $\iota(L') \in \Sigma$. The map $\iota: \Sigma' \to \Sigma$ induces a map of abstract simplicial complexes $\sigma: \Delta' \to \Delta$. Any map $\sigma: \Delta' \to \Delta$ which arises in this way is called a \emph{geometric simplicial subdivision} of $\Delta$. The crucial properties of such a map motivate the (semi-)abstract notion of simplicial subdivision
\cite[Section~2]{Sta92} which appears in the following definition. The notion of homology subdivision (over a fixed field $\kk$), adopted in \cite{Ath12a}, is a technical extension of that of topological subdivision which is often convenient (for instance,
for the setup of Section~\ref{sec:flag}) and sometimes essential
(for instance, in Theorem~\ref{thm:crosspolytope}). We will denote
by $2^V$ the abstract simplex consisting of all subsets of a set
$V$.

\begin{definition} \label{def:sub}
Let $\Delta$ be a simplicial complex. A (finite, simplicial) \emph{homology subdivision} of $\Delta$ (over $\kk$) is a simplicial complex $\Delta'$ together with a map $\sigma: \Delta' \to \Delta$ such that the following hold for every $F \in \Delta$: {\rm (a)} the set $\Delta'_F := \sigma^{-1} (2^F)$ is a subcomplex of $\Delta'$ which is a homology ball (over $\kk$) of dimension $\dim(F)$; and {\rm (b)} $\sigma^{-1} (F)$ consists of the interior faces of $\Delta'_F$.

Such a map $\sigma$ is said to be a \emph{topological subdivision} if $\Delta'_F$ is homeomorphic to a ball of dimension $\dim(F)$ for every $F \in \Delta$.
\end{definition}

\begin{remark}
To justify adopting this level of generality, let us mention that it is often more convenient to work in the setting of
Definition~\ref{def:sub}, rather than with the traditional notion of geometric subdivision, for the same reasons that abstract
simplicial complexes are often preferable to their geometric counterparts.
\end{remark}

Given faces $E \in \Delta'$ and $F \in \Delta$, the face $\sigma(E)$
of $\Delta$ is called the \emph{carrier} of $E$; the simplicial complex $\Delta'_F$ is called the \emph{restriction} of $\Delta'$
to $F$. The restriction of a homology (respectively, topological) subdivision $\sigma: \Delta' \to \Delta$ to $\Delta'_F$ is a homology (respectively, topological) subdivision of the simplex $2^F$ for every $F \in \Delta$. Recall that a subcomplex $\Gamma$ of a simplicial complex $\Delta$ is \emph{vertex-induced} if $\Gamma$ contains every face of $\Delta$ whose vertices belong to $\Gamma$.

\begin{definition} \label{def:quasi}
Let $\Delta', \Delta$ be simplicial complexes and $\sigma: \Delta' \to \Delta$ be a homology subdivision.
\begin{itemize}
\itemsep=0pt
\item[{\rm (a)}]
{\rm (\cite[Definition 2.4]{Ath12a})} The subdivision $\sigma$ is
\emph{vertex-induced} if $\Delta'_F$ is a vertex-induced subcomplex of $\Delta'$ for every $F \in \Delta$.
\item[{\rm (b)}]
{\rm (\cite[Definition 4.1]{Sta92})} The subdivision $\sigma$ is \emph{quasigeometric} if there do not exist $E \in \Delta'$ and face $F \in \Delta$ of dimension smaller than $\dim(E)$, such that the carrier of every vertex of $E$ is contained in $F$.
\end{itemize}
\end{definition}

For homology (respectively, topological) subdivisions we have the hierarchy of properties: geometric $\Rightarrow$ vertex-induced $\Rightarrow$ quasigeometric $\Rightarrow$ homology (respectively, topological). An example discussed on \cite[p.~468]{Cha94} shows that the first implication is strict. The following example shows that the other implications are strict as well.

  \begin{figure}
  \includegraphics[width=5in]{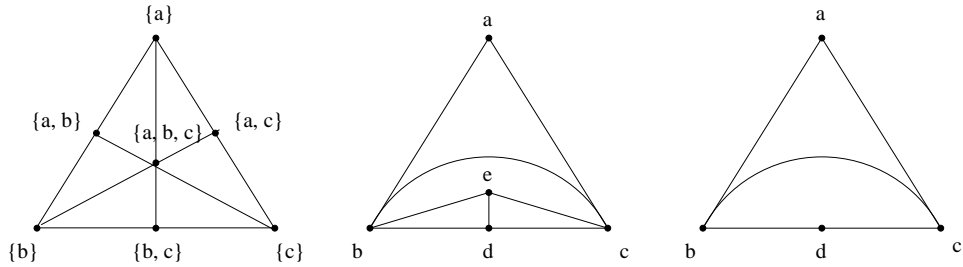}
  \caption{Three topological subdivisions of the 2-simplex}
  \label{fig:2dsub}
  \end{figure}

\begin{example} \label{ex:2dsub}
Let $V = \{a, b, c\}$ be a three element set. Figure~\ref{fig:2dsub} shows three topological subdivisions of the simplex
$2^V$. The leftmost (barycentric subdivision of $2^V$) is geometric, the middle one is quasigeometric but not vertex-induced and the rightmost is not quasigeometric. To view the latter, for instance, in the framework of Definition~\ref{def:sub}, one should take as $\Delta'$ the simplicial complex consisting of all subsets of $\{a, b, c\}$ and $\{b, c, d\}$ and set $\Delta = 2^V$. Then the carrier of all three faces of $\Delta'$ containing $\{b, c\}$ should be defined as $V$, the carrier of $\{b, d\}$, $\{d\}$ and $\{c, d\}$ should be defined as $\{b, c\}$ and all other faces of $\Delta'$ should be carried to themselves.
\end{example}

\subsection{The local $h$-vector}
\label{subsec:localh}

The local $h$-vector of a topological subdivision of a simplex is introduced in \cite[Definition~2.1]{Sta92}; it refines an invariant (combinatorial Newton number) considered in \cite{GZK91}. We record the definition for homology subdivisions.

\begin{definition}
\label{def:localh}
Let $V$ be a set with $d$ elements and let $\Gamma$ be a homology subdivision of the simplex $2^V$. The polynomial $\ell_V (\Gamma, x) = \ell_0 + \ell_1 x + \cdots + \ell_d x^d$ defined by
\begin{equation} \label{eq:deflocalh}
\ell_V (\Gamma, x) \ = \sum_{F \subseteq V} (-1)^{d - |F|} \,
h (\Gamma_F, x)
\end{equation}
is the \emph{local $h$-polynomial} of $\Gamma$ (with respect to $V$). The sequence $\ell_V (\Gamma) = (\ell_0,\dots,\ell_d)$ is the \emph{local $h$-vector} of $\Gamma$ (with respect to $V$).
\end{definition}

The following theorem, stated for homology subdivisions, shows that local $h$-vectors have surprisingly nice properties (see
Theorems~3.2, 3.3 and 5.2 and Corollary~4.7 in \cite{Sta92}). As usual, we denote by $\link_\Delta (F)$ the link of a simplicial complex $\Delta$ at the face $F \in \Delta$.

\begin{theorem} \label{thm:stalocal} {\rm (\cite{Sta92})}
\begin{itemize}
\item[{\rm (a)}]
For every homology subdivision $\Delta'$ of a pure simplicial complex $\Delta$,

\vspace*{-0.1in}
\begin{eqnarray*}
h (\Delta', x) &=& \sum_{F \in \Delta} \ell_F (\Delta'_F, x) \, h (\link_\Delta (F), x) \label{eq:hformula1} \\
    & & \nonumber \\
    &=& h (\Delta, x) \ \ + \ \sum_{F \in \Delta, \, |F| \ge 2} \ell_F (\Delta'_F, x) \, h (\link_\Delta (F), x). \label{eq:hformula2}
\end{eqnarray*}
\item[{\rm (b)}]
The local $h$-vector $\ell_V (\Gamma)$ is symmetric for every homology subdivision $\Gamma$ of the simplex $2^V$, i.e., we have $\ell_i = \ell_{d-i}$ for $0 \le i \le d$.
\item[{\rm (c)}]
The local $h$-vector $\ell_V (\Gamma)$ has nonnegative coordinates for every quasigeometric homology subdivision $\Gamma$ of the simplex
$2^V$.
\item[{\rm (d)}]
The local $h$-vector $\ell_V (\Gamma)$ is unimodal for every regular subdivision $\Gamma$ of the simplex $2^V$, i.e., we have $\ell_0 \le \ell_1 \le \cdots \le \ell_{\lfloor d/2 \rfloor}$.
\end{itemize}
\end{theorem}

\begin{remark} \label{rem:counterex}
Part (h) of \cite[Example~3.2]{Sta92} (due to Clara Chan) and
\cite[Example~3.4]{Ath12a} describe a non-quasigeometric subdivision and a quasigeometric but not vertex-induced subdivision, respectively, of the 3-simplex into two and five topological 3-simplices. These subdivisions are 3-dimensional analogues of the rightmost and middle subdivisions in Figure~\ref{fig:2dsub} and have local $h$-polynomials $-x^2$ and $x + x^3$, respectively. Thus, these examples show that the assumption in Theorem~\ref{thm:stalocal} (c) that $\Gamma$ is quasigeometric is essential and that there exist quasigeometric subdivisions with non-unimodal local $h$-vectors, contrary to what is suggested in \cite[Conjecture~5.4]{Sta92} \cite[Conjecture~5.11]{Sta94} \cite[Section~6]{Cha94} \cite[p.~134]{StaCCA}.
\end{remark}

\begin{example} \label{ex:2^n-1facets}
Theorem~\ref{thm:stalocal} is strong enough to determine the local
$h$-vector in certain situations. Consider, for instance, a quasigeometric homology subdivision $\Gamma$ of a $(d-1)$-dimensional simplex $2^V$ and assume that the restriction $\Gamma_F$ has exactly $2^{\dim(F)}$ maximal faces for every nonempty face $F$ of $2^V$. Recalling that the number of maximal faces of a simplicial complex $\Delta$ is equal to $h(\Delta, 1)$ and setting $x=1$ in the defining equation (\ref{eq:deflocalh}), an easy computation shows that $\ell_V (\Gamma, 1)$ is equal to 1, if $d$ is even, and to 0 otherwise. Parts (b) and (c) of Theorem~\ref{thm:stalocal} then imply that
  $$ \ell_V (\Gamma, x) \ = \ \begin{cases} x^{d/2}, & \text{if $d$ is even}, \\
      0, & \text{if $d$ is odd}. \end{cases} $$
\end{example}

The following statement is a consequence of parts (a) and (c) of Theorem~\ref{thm:stalocal} and partially answers the question of Kalai and Stanley on the monotonicity of $h(\Delta)$, mentioned in the introduction.

\begin{corollary} \label{cor:hmonotone}
{\rm (\cite[Theorem~4.10]{Sta92})}
We have $h(\Delta') \ge h(\Delta)$ coordinatewise for every
Cohen-Macaulay (or even Buchsbaum) simplicial complex $\Delta$ and every quasigeometric homology subdivision $\Delta'$ of $\Delta$.
\end{corollary}

We end this subsection with some major open problems from \cite{Ath12a, Cha94, Sta92}.

\begin{question} \label{que:unimodal}
{\rm (\cite[Question~3.5]{Ath12a})} Is the local $h$-vector $\ell_V (\Gamma)$ unimodal for every geometric
(or even vertex-induced homology) subdivision $\Gamma$ of the simplex $2^V$?
\end{question}

\begin{problem} \label{pro:characterize}
Characterize local $h$-vectors of quasigeometric and of
vertex-induced homology subdivisions of the simplex (see \cite{Cha94} for related results and Remark~\ref{rem:counterex}).
\end{problem}

\begin{problem} \label{pro:zero}
{\rm (\cite[Problem~4.13]{Sta92})}
Characterize explicitly the quasigeometric homology subdivisions $\Gamma$ of the simplex $2^V$ for which $\ell_V (\Gamma, x) = 0$.
\end{problem}

\begin{question} \label{que:topo}
{\rm (\cite[Conjecture~4.11]{Sta92})}
Can the hypothesis in Corollary~\ref{cor:hmonotone} that $\Delta'$ is quasigeometric be dropped?
\end{question}

\subsection{The relative local $h$-polynomial}
\label{subsec:rlocalh}

The relative local $h$-polynomial was introduced in \cite[Remark~3.7]{Ath12a} for homology subdivisions of a simplex and independently (in a different level of generality) in \cite[Section~3]{Ni12} for regular triangulations of polytopes.

\begin{definition} \label{def:relative}
Let $\sigma: \Gamma \to 2^V$ be a homology subdivision of a
$(d-1)$-dimensional simplex $2^V$ and let $E \in \Gamma$. The polynomial
\begin{equation} \label{eq:deflocalhrel}
\ell_V (\Gamma, E, x) \ = \sum_{\sigma(E) \subseteq F \subseteq V}
(-1)^{d - |F|} \, h (\link_{\Gamma_F} (E), x)
\end{equation}
is the \emph{relative local $h$-polynomial} of $\Gamma$ (with respect to $V$) at $E$.
\end{definition}

The polynomial $\ell_V (\Gamma, E, x)$ reduces to $\ell_V (\Gamma,
x)$ for $E = \varnothing$. The following statement shows that $\ell_V (\Gamma, E, x)$ enjoys properties analogous to those of $\ell_V (\Gamma, x)$.

\begin{theorem} \label{thm:relative}
{\rm (\cite[Section~3]{Ath12a})}
Let $V$ be a $d$-element set.
\begin{itemize}
\item[{\rm (a)}]
For every homology subdivision $\Gamma$ of the simplex $2^V$ and every homology subdivision $\Gamma'$ of $\Gamma$,

\vspace*{-0.1in}
\begin{eqnarray*}
\ell_V (\Gamma', x) &=&
\sum_{E \in \Gamma} \ell_E (\Gamma'_E, x) \, \ell_V (\Gamma, E, x) \\
&=& \ell_V (\Gamma, x) \ \ + \ \sum_{E \in \Gamma, \, |E| \ge 2} \ell_E (\Gamma'_E, x) \, \ell_V (\Gamma, E, x).
\end{eqnarray*}

\item[{\rm (b)}]
The polynomial $\ell_V (\Gamma, E, x)$ has symmetric coefficients, meaning we have
$x^{d-|E|} \, \ell_V (\Gamma, E, 1/x) = \ell_V (\Gamma, E, x)$,
for every homology subdivision $\Gamma$ of the simplex $2^V$ and every $E \in \Gamma$.
\item[{\rm (c)}]
The polynomial $\ell_V (\Gamma, E, x)$ has nonnegative coefficients for every quasigeometric homology subdivision $\Gamma$ of the simplex $2^V$ and every $E \in \Gamma$.
\end{itemize}
\end{theorem}

\noindent
\emph{Proof (sketch)}. Part (a) was proven in \cite[Proposition~3.6]{Ath12a}. We comment on the proofs of parts (b) and (c), which were omitted in \cite{Ath12a}. For part (b) one can adapt the proof of \cite[Theorem~4.2]{Ath12b}; a detailed description of the argument is recorded in \cite[Section~3.5]{Sav13}. The special case $E =
\varnothing$ of part (c) is equivalent to Theorem~\ref{thm:stalocal} (c). The general case follows by the argument in the proof of
\cite[Theorem~5.1]{Ath12b} (generalizing that in the proof of
\cite[Theorem~4.6]{Sta92}), where the role of $\Delta$ in that proof is played by $\link_\Gamma (E)$ (which is a homology sphere, if $\sigma(E) = V$, and a homology ball otherwise), the role of $d$ is played by $d - |E| = \dim \link_\Gamma (E) + 1$ and the role of $e$ is played by the rank $d - |\sigma(E)|$ of the interval $[\sigma(E), V]$ in the lattice of subsets of $V$.
\qed

\medskip
A combinatorial application of part (a) of Theorem~\ref{thm:relative} is given in Example~\ref{ex:relativeapp}.

\begin{question} \label{que:relunimodal}
Does the relative local $h$-polynomial $\ell_V (\Gamma, E, x)$ have  unimodal coefficients for every regular (or just geometric)
subdivision $\Gamma$ of the simplex $2^V$ and every $E \in \Gamma$?
\end{question}

The following monotonicity property of local $h$-vectors is a consequence of Theorem~\ref{thm:stalocal} (c) and parts (a) and (c) of Theorem~\ref{thm:relative}.

\begin{corollary} {\rm (\cite[Remark~3.7]{Ath12a})}
We have $\ell_V (\Gamma', x) \ge \ell_V (\Gamma, x)$ coefficientwise for every quasigeometric homology subdivision $\Gamma$ of the simplex $2^V$ and every quasigeometric homology subdivision $\Gamma'$ of $\Gamma$.
\end{corollary}

\section{Flag homology spheres and $\gamma$-vectors}
\label{sec:flag}

A simplicial complex $\Delta$ is called \emph{flag} if every minimal nonface of $\Delta$ has two elements; see \cite[Section~III.4]{StaCCA} for a discussion of this fascinating class of complexes. This section surveys some recent developments, extending part of the  theory of \cite{Sta92} to subdivisions of flag homology spheres and their $\gamma$-vectors. The following notion was introduced in
\cite{Ath12a}.
\begin{definition} \label{def:flagsub}
{\rm (\cite[Definition~2.4]{Ath12a})} We say that a homology subdivision $\Delta'$ of a simplicial complex $\Delta$ is \emph{flag} if the restriction $\Delta'_F$ is a flag complex for every face $F \in \Delta$.
\end{definition}

We recall (see, for instance, Section~II.6) that the $h$-polynomial of a homology sphere (more generally, of an Eulerian simplicial complex) $\Delta$ has symmetric coefficients, i.e., we have $h_i(\Delta) = h_{d-i}(\Delta)$ for $0 \le i \le d$, where $\dim(\Delta) = d-1$. Thus, it can be written uniquely in the form
\begin{equation} \label{eq:defgamma}
h (\Delta, x) \  = \, \sum_{i=0}^{\lfloor d/2 \rfloor} \gamma_i x^i
(1+x)^{d-2i}.
\end{equation}
The sequence $\gamma (\Delta) = (\gamma_0,
\gamma_1,\dots,\gamma_{\lfloor d/2 \rfloor})$ is the
\emph{$\gamma$-vector} and $\gamma (\Delta, x) = \sum_{i=0}^{\lfloor d/2 \rfloor} \gamma_i x^i$ is the \emph{$\gamma$-polynomial} of
$\Delta$ \cite{Ga05}. The following conjecture of Gal, which extends an earlier conjecture of Charney and Davis~\cite{CD95}, can be viewed as a Generalized Lower Bound Conjecture for flag homology spheres; note that it implies the unimodality of $h(\Delta)$.
\begin{conjecture} \label{conj:gal}
{\rm (\cite[Conjecture~2.1.7]{Ga05})} The vector $\gamma (\Delta)$ has nonnegative coordinates for every flag homology sphere $\Delta$.
\end{conjecture}

Gal's conjecture is known to hold (essentially by a result of
\cite{DO01}) for spheres of dimension at most four \cite[Section~2]{Ga05}, but is open even for (flag) boundary complexes of simplicial polytopes in higher dimensions. The monotonicity question for
$\gamma$-vectors was raised first by Postnikov, Reiner and
Williams~\cite[Conjecture~14.2]{PRW08} for geometric subdivisions and later by this author for more general subdivisions as follows.
\begin{conjecture} \label{conj:gammamonotone}
{\rm (\cite[Conjecture~1.4]{Ath12a})} We have $\gamma (\Delta')
\ge \gamma (\Delta)$ coordinatewise for every flag homology sphere $\Delta$ and every flag vertex-induced homology subdivision $\Delta'$ of $\Delta$.
\end{conjecture}

To transfer the theory of \cite{Sta92} to subdivisions of flag homology spheres and $\gamma$-vectors, one needs the following concept.
\begin{definition} \label{def:localgamma}
{\rm (\cite[Definition~5.1]{Ath12a})} Let $V$ be a $d$-element set and $\Gamma$ be a homology subdivision of the simplex $2^V$. The polynomial $\xi_V (\Gamma, x) = \sum_{i=0}^{\lfloor d/2 \rfloor} \xi_i x^i$ uniquely defined by
\begin{equation} \label{eq:defxi}
\ell_V (\Gamma, x) \  = \, \sum_{i=0}^{\lfloor d/2 \rfloor} \xi_i x^i (1+x)^{d-2i}
\end{equation}
is the \emph{local $\gamma$-polynomial} of $\Gamma$ (with respect to $V$). The sequence $\xi_V (\Gamma) =
(\xi_0,\dots,\xi_{\lfloor d/2 \rfloor})$ is the \emph{local
$\gamma$-vector} of $\Gamma$ (with respect to $V$).
\end{definition}

\begin{theorem} \label{thm:athalocal}
{\rm (\cite[Section~5]{Ath12a})}
\begin{itemize}
\item[{\rm (a)}]
For every homology sphere (or even pure Eulerian simplicial complex) $\Delta$ and every homology subdivision $\Delta'$ of $\Delta$ we have

\vspace*{-0.1in}
\begin{eqnarray*}
\gamma (\Delta', x) &=& \sum_{F \in \Delta} \xi_F (\Delta'_F, x) \, \gamma (\link_\Delta (F), x) \\
    & & \\
    &=& \gamma (\Delta, x) \ \ + \ \sum_{F \in \Delta, \, |F| \ge 2} \xi_F (\Delta'_F, x) \, \gamma (\link_\Delta (F), x) .
\end{eqnarray*}
\item[{\rm (b)}]
The vector $\xi_V (\Gamma)$ has nonnegative coordinates for every flag vertex-induced homology subdivision $\Gamma$ of a simplex $2^V$ of dimension at most three and for every subdivision $\Gamma$ which can be obtained from the trivial subdivision of a simplex $2^V$ of any dimension by successive stellar subdivisions on edges.
\item[{\rm (c)}]
Conjecture~\ref{conj:gammamonotone} holds whenever $\Delta$ (and hence $\Delta'$) has dimension at most four.
\end{itemize}
\end{theorem}

These results motivated the following conjecture.
\begin{conjecture} \label{conj:main}
{\rm (\cite[Conjecture~5.4]{Ath12a})} The local $\gamma$-vector $\xi_V (\Gamma)$ has nonnegative coordinates for every flag
vertex-induced homology subdivision $\Gamma$ of the simplex $2^V$.
\end{conjecture}

As Theorem~\ref{thm:athalocal} (a) shows, Conjecture~\ref{conj:main}
together with Gal's conjecture imply
Conjecture~\ref{conj:gammamonotone}. The following statement, which is of independent interest, shows that in fact we have the implications Conjecture~\ref{conj:main} $\Rightarrow$ Conjecture~\ref{conj:gammamonotone} $\Rightarrow$ Conjecture~\ref{conj:gal}. We denote by $\Sigma_{d-1}$ the boundary complex of the $d$-dimensional cross-polytope (this complex plays for flag homology spheres the role played by the boundary complex of the simplex for all homology spheres; it satisfies $\gamma(\Sigma_{d-1}, x) = 1$).
\begin{theorem} \label{thm:crosspolytope}
{\rm (\cite[Theorem~1.5, Corollary~5.5]{Ath12a})}
 Every flag $(d-1)$-dimensional homology sphere $\Delta$ is a
vertex-induced (hence quasigeometric and flag) homology subdivision of $\Sigma_{d-1}$. Moreover, we have
  \begin{equation} \label{eq:gamma-xi}
    \gamma (\Delta, x) \ = \ \sum_{F \in \Sigma_{d-1}} \, \xi_F (\Delta_F, x).
  \end{equation}
\end{theorem}

%%We end this section with the following open problem.

\begin{problem} \label{pro:mainconjecturebary} \rm
(\cite[Question~6.2]{Ath12a}) Study Conjecture~\ref{conj:main} for barycentric subdivisions of geometric (simplicial or, more generally, polyhedral) subdivisions of the simplex.
\end{problem}

\section{Examples}
\label{sec:ex}

This section surveys recent results which give combinatorial interpretations, or other formulas, for the local $h$-vector and local $\gamma$-vector for certain interesting families of flag geometric subdivisions (including various barycentric, edgewise and cluster subdivisions) of the simplex and thus verify Conjecture~\ref{conj:main} in these cases. One previously unpublished result
(Theorem~\ref{thm:edgewise}) is included. Throughout this section, $V$ stands for an $n$-element set.

The (first, simplicial) barycentric subdivision of the simplex $2^V$, denoted here by $\sd(2^V)$, is a flag subdivision of well known importance in algebraic topology and combinatorics (the $n=3$ case appears in the left part of Figure~\ref{fig:2dsub}). The interpretation to the local $h$-polynomial of $\sd(2^V)$, given in (\ref{eq:localbary}), is due to Stanley~\cite[Section~2]{Sta92}.
The interpretation to the local $\gamma$-polynomial, given in
(\ref{eq:localgammabary}), has appeared (implicitly or explicitly) in various contexts in independent works by several authors; methods used include the Foata-Sch\"utzenberger valley hopping \cite[Section~4]{AS12}, symmetric function expansions \cite[Section~4]{LSW12} \cite[Theorem~7.3]{SW10}, Rees product homology \cite{LSW12} and continued fractions \cite{SZ12}. These references provide several interesting generalizations and refinements. We denote by $\mathfrak{S}_n$ the symmetric group of permutations of $\{1,
2,\dots,n\}$ and by $\dD_n$ the set of derangements (permutations without fixed points) in $\mathfrak{S}_n$. We recall that an \emph{excedance} of $w \in \mathfrak{S}_n$ is an index $1 \le i \le n$ such that $w(i) > i$. A \emph{descending run} of $w \in \mathfrak{S}_n$ is a maximal string $\{a, a+1,\dots,b\}$ of integers, such that $w(a) > w(a+1) > \cdots > w(b)$.

\begin{theorem} \label{thm:bary} {\rm (\cite[Proposition~2.4]{Sta92} \cite[Theorem~1.4]{AS12} \cite[Equation~(1.3) and Theorem~3.3]{LSW12} \cite[Remark~5.5]{SW10} \cite[Corollary~9, Theorem~11]{SZ12})}
The local $h$-polynomial of the barycentric
subdivision of the $(n-1)$-dimensional simplex $2^V$ can be expressed as
\begin{eqnarray} \label{eq:localbary}
\ell_V (\sd(2^V), x) &=& \sum_{w \in \dD_n} x^{\ex(w)} \\
& & \nonumber \\
&=& \label{eq:localgammabary}
\sum_{i=0}^{\lfloor n/2 \rfloor} \xi_{n, i} \, x^i (1+x)^{n-2i},
\end{eqnarray}
where $\ex(w)$ is the number of excedances of $w \in \mathfrak{S}_n$ and $\xi_{n,i}$ stands for the number of permutations $w \in \mathfrak{S}_n$ with $i$ descending runs and no descending run of size one.
\end{theorem}

\begin{remark} \label{rem:dn(x)}
The right-hand side of (\ref{eq:localbary}) is usually referred to as the \emph{derangement polynomial} of order $n$ and is denoted by $d_n (x)$. It was first considered and shown to be unimodal by Brenti~\cite{Bre90}. Theorem~\ref{thm:bary} gives a combinatorial proof of the unimodality of this polynomial. Another such proof was given by Stembridge~\cite{Ste92}.
\end{remark}

  \begin{figure}
  \includegraphics[width=4.3in]{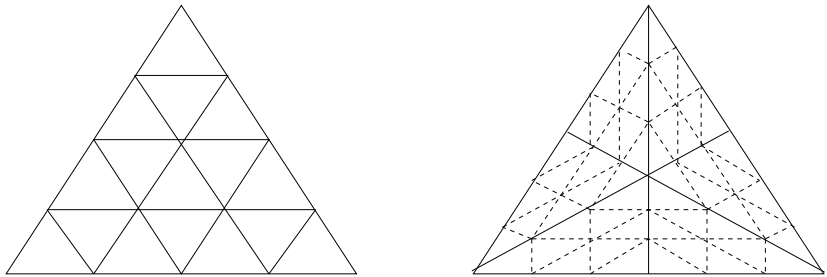}
  \caption{Two edgewise subdivisions}
  \label{fig:edge34}
  \end{figure}

Given a positive integer $r$, the $r$th edgewise subdivision $\Delta^{\langle r \rangle}$ is another standard way to subdivide a simplicial complex $\Delta$ by which each face $F \in \Delta$ is subdivided into $r^{\dim(F)}$ faces of the same dimension. We refer to \cite{BrW09, BR05} for the precise definition and for a discussion of the long history of this construction in mathematics. Figure~\ref{fig:edge34} shows the fourth edgewise subdivision, and the third edgewise subdivision of the barycentric subdivision, of the
2-dimensional simplex. Theorem~\ref{thm:bary} can be generalized as follows. Consider the group $\ZZ_r \wr \mathfrak{S}_n$ of $r$-colored permutations, meaning permutations in $\mathfrak{S}_n$ with each coordinate in their one--line notation colored with an element of the set $\{0, 1,\dots,r-1\}$. The \emph{flag excedance}~\cite{BG06} of an element $w \in \ZZ_r \wr \mathfrak{S}_n$ is defined as $\fex(w) = r \cdot \ex_A (w) + \csum(w)$, where $\ex_A (w)$ is the number of coordinates $i \in \{1, 2,\dots,n\}$ of $w$ of zero color such that $w(i) > i$ and $\csum(w)$ stands for the sum of the colors of all coordinates of $w$. We denote by $(\dD^r_n)^\circ$ the set of derangements (elements without fixed points of zero color) $w \in \ZZ_r \wr \mathfrak{S}_n$ for which $\csum(w)$ is divisible by $r$. Decreasing runs for elements of $\ZZ_r \wr \mathfrak{S}_n$ are defined as for those of $\mathfrak{S}_n$, after the $r$-colored integers have been totally ordered in a standard lexicographic way
(see~\cite{Ath13} for more explanation).

\begin{theorem} \label{thm:edgebary}
{\rm (\cite[Theorem~1.2]{Ath13})}
The local $h$-polynomial of the $r$th edgewise subdivision of the barycentric subdivision of the $(n-1)$-dimensional simplex $2^V$ can be expressed as
\begin{eqnarray} \label{eq:localedgebary}
\ell_V (\sd(2^V)^{\langle r \rangle}, x) &=& \sum_{w \in
(\dD^r_n)^\circ} x^{\fex (w) / r} \\
& & \nonumber \\
&=& \label{eq:localgammaedgebary}
\sum_{i=0}^{\lfloor n/2 \rfloor} \xi^+_{n, r, i} \, x^i (1+x)^{n-2i},
\end{eqnarray}
where $\xi^+_{n, r, i}$ is the number of elements of $\ZZ_r \wr \mathfrak{S}_n$ with $i$ descending runs, no descending run of size one and last coordinate of zero color.
\end{theorem}

\begin{remark} \label{rem:reesproof}
The right-hand side of (\ref{eq:localedgebary}) is closely related to the derangement polynomial for the group $\ZZ_r \wr \mathfrak{S}_n$, introduced and studied by Chow and Mansour~\cite{CM10}; see
\cite[Theorem~1.3]{Ath13}. The proof of (\ref{eq:localgammaedgebary}) given in~\cite{Ath13} uses a result of \cite{LSW12} on Rees product homology; a more direct combinatorial proof should be possible.
\end{remark}

\begin{figure}
  \includegraphics[width=4.3in]{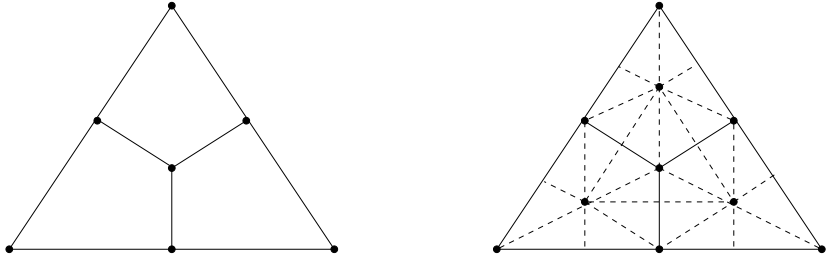}
  \caption{The cubical barycentric subdivision of the 2-simplex and
  its barycentric subdivision}
  \label{fig:K3}
  \end{figure}

\begin{remark} \label{rem:cubical}
The formulas of the special case $r=2$ of Theorem~\ref{thm:edgebary} are also valid for the local $h$-polynomial of the (simplicial) barycentric subdivision of the cubical barycentric subdivision of the simplex $2^V$ (studied in \cite[Chapter~3]{Sav13} and shown in Figure~\ref{fig:K3} for $n=3$; this example motivated
Theorem~\ref{thm:edgebary}). Indeed, let us denote this subdivision by $K_n$. One can show (see \cite[Section~3.6]{Sav13}) that $K_n$ is a subdivision of $\sd(2^V)$ and that the restriction of $K_n$ to each nonempty face $F \in \sd(2^V)$ has exactly $2^{\dim(F)}$ maximal faces. Since the same is true for $\sd(2^V)^{\langle 2 \rangle}$, as already discussed, it follows from Theorem~\ref{thm:relative} (a) and Example~\ref{ex:2^n-1facets} that $\ell_V (K_n, x) = \ell_V
(\sd(2^V)^{\langle 2 \rangle}, x)$.
\end{remark}

We now turn attention to the $r$th edgewise subdivision of the simplex $2^V$. We denote by $\sS\wW (n, r)$ the set of all words
$w = (w_0, w_1,\dots,w_n) \in \{0, 1,\dots,r-1\}^{n+1}$ satisfying $w_0 = w_n = 0$ and having no two consecutive entries equal.
We call an index $i \in \{0, 1,\dots,n-1\}$ an \emph{ascent} of the word $w$ if $w_i < w_{i+1}$ and an index $i \in \{1, 2,\dots,n-1\}$ a
\emph{double ascent} of $w$ if $w_{i-1} < w_i < w_{i+1}$ (descents and double descents are similarly defined). The elements of $\sS\wW (n,r)$ are the Smirnov words in the set $\{0, 1,\dots,r-1\}^{n+1}$ with zero first and last coordinate; see \cite[Section~5]{LSW12}, where some results of similar flavor with the following theorem appear. We will also denote by ${\rm E}_r$ the linear operator on $\RR[x]$ defined by setting ${\rm E}_r (x^m) = x^{m/r}$, if $m$ is divisible by $r$, and ${\rm E}_r (x^m) = 0$ otherwise. Note that the special case $r=2$ of the following theorem agrees with Example~\ref{ex:2^n-1facets}.

\begin{theorem} \label{thm:edgewise}
The local $h$-polynomial of the $r$th edgewise subdivision of the
$(n-1)$-dimensional simplex $2^V$ can be expressed as
\begin{eqnarray} \label{eq:localedge}
\ell_V ((2^V)^{\langle r \rangle}, x) &=& {\rm E}_r \, (x + x^2 + \cdots + x^{r-1})^n \ \ = \ \sum_{w \in \sS\wW(n, r)} x^{\asc (w)} \\
&=& \label{eq:localgammaedge}
\sum_{i=0}^{\lfloor n/2 \rfloor} \xi_{n, r, i} \, x^i (1+x)^{n-2i},
\end{eqnarray}
where $\asc (w)$ is the number of ascents of $w \in \sS\wW(n, r)$ and $\xi_{n, r, i}$ stands for the number of words $w = (w_0,
w_1,\dots,w_n) \in \sS\wW(n, r)$ with exactly $i$ ascents which have the following property: for every double ascent $j$ of $w$ there exists a double descent $k > j$ such that $w_j = w_k$ and $w_j \le w_r$ for all $j < r < k$.
\end{theorem}

\begin{remark} \label{rem:edgegammatransform}
The transformation of the $\gamma$-vector of a homology sphere after barycentric subdivision was studied in~\cite{NPT11}. It would be interesting to study the corresponding transformation for the $r$th edgewise subdivision.
\end{remark}

\begin{example} \label{ex:relativeapp}
As an application of Theorem~\ref{thm:relative} (a), we give another formula for the local $h$-polynomial considered in
Theorem~\ref{thm:edgebary} as follows. Let $\Gamma$ be the barycentric subdivision of the simplex $2^V$ and $E = \{S_1,
S_2,\dots,S_k\}$ be a face of $\Gamma$, where $S_1 \subset S_2 \subset \cdots \subset S_k \subseteq V$ are nonempty sets. Using the defining equation (\ref{eq:deflocalhrel}), one can show (see
\cite[Example~3.5.2]{Sav13}) that
  \begin{equation} \label{eq:baryrel}
    \ell_V (\sd(2^V), E, x) \ = \ d_{n_0}(x) \, A_{n_1}(x) A_{n_2}(x)
    \cdots A_{n_k}(x),
  \end{equation}
where $d_{n_0}(x)$ is the derangement polynomial discussed in
Remark~\ref{rem:dn(x)}, $n_0 = |V \sm S_k|$ and $n_i = |S_i \sm
S_{i-1}|$ for $1 \le i \le k$ (with the convention $S_0 =
\varnothing$). Theorem~\ref{thm:relative} (a) applied to $\Gamma = \sd(2^V)$ and $\Gamma' = \sd(2^V)^{\langle r \rangle}$ and
Equations~(\ref{eq:localedge}) and (\ref{eq:baryrel}) imply that
\[ \ell_V (\sd(2^V)^{\langle r \rangle}, x) \ = \
\sum {n \choose n_0, n_1,\dots,n_k} \, {\rm E}_r \, (x + \cdots + x^{r-1})^k \, d_{n_0} (x) \, A_{n_1} (x) \cdots A_{n_k}(x), \]
where the sum ranges over all $k \in \NN$ and over all sequences
$(n_0, n_1,\dots,n_k)$ of integers which satisfy $n_0 \ge 0$, $n_1,\dots,n_k \ge 1$ and sum to $n$. For the special case $r=2$, more formulas, recurrences and generating functions appear in
\cite[Chapter~3]{Sav13}.
\end{example}

Our final example comes from the theory of cluster algebras. Let $\Phi$ be a finite root system with positive system $\Phi^+$ and simple system $\Pi$. The cluster complex $\Delta(\Phi)$ \cite{FZ03} is a flag simplicial sphere on the vertex set $\Phi^+ \cup (-\Pi)$
with remarkable enumerative properties. The restriction of $\Delta(\Phi)$ to the vertex set $\Phi^+$ is naturally a flag geometric subdivision $\Gamma(\Phi)$ of the simplex $2^\Pi$, termed in
\cite{AS12} as the \emph{cluster subdivision}. Since the enumerative combinatorics of $\Delta(\Phi)$ has proved to be so interesting
(see \cite{FZ03} and references in \cite{AS12}), one may expect that the local $h$-vector and local $\gamma$-vector of $\Gamma
(\Phi)$ are interesting as well. We record the relevant result for the root system of type $A_n$; similar interpretations or computations for the other root systems appear in \cite{AS12}.

\begin{theorem} \label{thm:cluster} {\rm (\cite{AS12})} The local
$h$-polynomial of the cluster subdivision of Cartan-Killing type
$A_n$ can be expressed as
\begin{eqnarray} \label{eq:localcluster}
\ell_\Pi (\Gamma (\Phi), x) &=& \sum_{i=0}^n \ell_i (\Phi) x^i \\
& & \nonumber \\
&=& \label{eq:localgammacluster}
\sum_{i=1}^{\lfloor n/2 \rfloor} \frac{1}{n-i+1} {n \choose i}
{n-i-1 \choose i-1} \, x^i (1+x)^{n-2i},
\end{eqnarray}
where $\ell_i (\Phi)$ is the number of noncrossing partitions $\pi$ of $\{1, 2,\dots,n\}$ with $i$ blocks, such that for every singleton block $\{b\}$ of $\pi$ there exist elements $a, c$ of some block of $\pi$ satisfying $a < b < c$.
\end{theorem}

We conclude with a few open problems. A 5-dimensional flag simplicial (polytopal) sphere whose $h$-polynomial is not real-rooted was constructed by Gal~\cite{Ga05}, thus disproving a stronger version of Conjecture~\ref{conj:gal}. Based on this example, one can prove the existence of a flag geometric (simplicial) subdivision of a simplex of dimension 7 whose local $h$-polynomial is not real-rooted.

\begin{problem} \label{pro:nonrealroot}
Find a flag geometric subdivision of a simplex of minimum possible dimension whose local $h$-polynomial is not real-rooted.
\end{problem}

\begin{question} \label{que:realroot}
Are the local $h$-polynomials of the subdivisions discussed in this section real-rooted?
\end{question}

\begin{problem} \label{pro:2bary}
Find a combinatorial interpretation for the local $h$-polynomial and the local $\gamma$-polynomial of the second barycentric subdivision of the simplex.
\end{problem}

\noindent
\emph{Note added in revision.} An affirmative answer to
Question~\ref{que:relunimodal} has been provided for regular subdivisions by E.~Katz and A.~Stapledon~\cite[Remark~6.5]{KaSt14}. An affirmative answer to Question~\ref{que:realroot} follows from the work of X.~Zhang~\cite{Zh95} for the barycentric subdivision, and has been claimed by P.~Zhang~\cite{Zh14} for the $r$th edgewise and cluster subdivisions, of the simplex.

\bibliographystyle{amsalpha}

\end{document}